\documentclass{amsart}

\usepackage{amssymb}

\usepackage{graphicx}

\usepackage[cmtip,all]{xy}

\usepackage{}
\usepackage{amsmath}
\usepackage{amsfonts}
\usepackage{amssymb}

\usepackage{mathrsfs}

\usepackage{amssymb,amsmath,amsxtra,graphicx}
\usepackage{times}
\usepackage{bbold}
\usepackage{stmaryrd}
\usepackage{array}
\usepackage{dsfont}
\usepackage{url}
\usepackage{tikz}
\usepackage[T1]{fontenc}

\newtheorem{theorem}{Theorem}
\newtheorem*{lem}{Lemma 0}

\theoremstyle{plain}

\newtheorem{corollary}[theorem]{Corollary}

\newtheorem{lemma}[theorem]{Lemma}

\begin{document}
\title{Two Pointwise Characterizations of the Peano Derivative}
\author{J. Marshall Ash}
\address{Department of Mathematics, DePaul University\\Chicago, IL 60614}
\email{mash@depaul.edu}
\author{Stefan Catoiu}
\address{Department of Mathematics, DePaul University\\Chicago, IL 60614}
\email{scatoiu@depaul.edu}
\author{Hajrudin Fejzi\'{c}}
\address{Department of Mathematics, California State University, San Bernardino, CA 92407}
\email{hfejzic@csusb.edu}
\thanks{May 8, 2024. This paper is in final form and no
version of it will be submitted for publication elsewhere.}

\begin{abstract} We provide the first two examples of sets of generalized Riemann derivatives  of orders up to $n$, $n\geq 2$, whose simultaneous existence for all functions~$f$ at~$x$ is equivalent to the existence of the $n$-th Peano derivative $f_{(n)}(x)$. In this way, we begin to understand how the theory of Peano derivatives can be explained exclusively in terms of generalized Riemann derivatives, a bold new principle in generalized differentiation.

In 1936, J. Marcinkiewicz and A. Zygmund showed that the existence of $f_{(n)}(x)$ is equivalent to the existence of both $f_{(n-1)}(x)$ and the $n$th generalized Riemann derivative $\widetilde{D}_nf(x)$, based at $x,x+h,x+2h,x+2^2h,\ldots ,x+2^{n-1}h$. Our first characterization of $f_{(n)}(x)$ is that its existence is equivalent to the simultaneous existence of $\widetilde{D}_1f(x),\ldots,\widetilde{D}_nf(x)$.
Our second characterization is that the existence of $f_{(n)}(x)$ is equivalent to the existence of $\widetilde{D}_1f(x)$ and of all $n(n-1)/2$ forward shifts,
\[
D_{k,j}f(x)=\lim_{h\rightarrow 0} h^{-k}\sum_{i=0}^k(-1)^i\binom ki f(x+(k+j-i)h),
\]
for $j=0,1,\ldots,k-2$, of the $k$-th Riemann derivatives $D_{k,0}f(x)$, for $k=2,\ldots ,n$.

The proof of the second result involves an interesting combinatorial algorithm that starts with consecutive forward shifts of an arithmetic progression and yields a geometric progression, using two set-operations: dilation and combinatorial Gaussian elimination.
This result proves a variant of a 1998 conjecture by Ginchev, Guerragio and Rocca, predicting the same outcome for backward shifts instead of forward shifts. The conjecture has been recently settled in \cite{AC2}, with a proof that has this variant's proof as a prerequisite.
\end{abstract}
\maketitle

\noindent 
\emph{Generalized Riemann derivatives}, or $\mathcal{A}$-derivatives, of order $n$ of a function $f$ at $x$ are given by limits of the form%
\[
D_{\mathcal{A}}f(x)=\lim_{h\rightarrow0}h^{-n}\sum_{i=0}^{m}A_{i}f\left(  x+a_{i}h\right)  ,
\]
where $m\geq n$ and the data vector $\mathcal{A}=\{A_0,\ldots ,A_m;a_0,\ldots ,a_m\}$ satisfies the $n$th \emph{Vandermonde conditions} $\sum_{i}A_{i}a_{i}^{j}=n!\delta_{jn}$, for $j=0,1,\dots,n$. If $a_i^j=0^0$, take this to be~1.
For example, the generalized Riemann derivative of order $n=1$ and defined by data vector $\mathcal{A}=\left\{  -1,1;0,1\right\}  $ is just ordinary differentiation, while the Schwarz derivative $S$ is the $\mathcal{A}$-derivative of order $n=2$ defined by $\mathcal{A}=\{1,-2,1;-1,0,1\}$. The $n$th forward \emph{Riemann derivative} $D_nf(x)$ is the $n$th generalized Riemann derivative $D_{\mathcal{A}}f(x)$, given by
\[\mathcal{A}=\left\{A_i=(-1)^i\binom ni;a_i=n-i\mid i=0,\ldots ,n\right\},\]
and the $n$th \emph{symmetric Riemann derivative} $D_n^sf(x)$ is the $n$th generalized Riemann derivative $D_{\mathcal{A}}f(x)$ defined by
$\mathcal{A}=\left\{A_i=(-1)^i\binom ni;a_i=\frac n2-i\,|\,i=0,\ldots ,n\right\}$. In particular, $D_1f(x)=f'(x)$ and $D_2^sf(x)=Sf(x)$.
All generalized Riemann derivatives considered in this paper have $m=n$, in which case, the non-zero outer coefficients $A_i$ are uniquely determined by the distinct inner coefficients $a_i$. Such derivatives are called \emph{exact}.

A function $f$ is $n$ times \emph{Peano differentiable} at~$x$, if it is Taylor aproximable to order~$n$ about~$x$, that is, if there exist numbers $f_{(0)}(x)=f(x),f_{(1)}(x),\ldots,f_{(n)}(x)$, called \emph{the first~$n$ Peano derivatives} of $f$ at $x$, such that
\[
f(x+h)=f_{(0)}(x)+f_{(1)}(x)h+\frac {f_{(2)}(x)}{2!}h^2+\cdots +\frac {f_{(n)}(x)}{n!}h^n+o(h^n),
\]
as $h$ approaches zero. The condition $f_{(0)}(x)=f(x)$ is to ensure that Peano differentiability of order~$0$ at~$x$ is the same as continuity at~$x$, hence, by Taylor's Theorem, that~$n$ times Peano differentiability is more general than~$n$ times ordinary differentiability for all~$f$ at~$x$. The above definition implies that $f_{(1)}(x)=f'(x)$. And in orders~$n$ at least two, the $n$th Peano derivative $f_{(n)}(x)$ is known to be strictly more general than the $n$th ordinary derivative $f^{(n)}(x)$.
Notice that the existence of the $n$th Peano derivative of~$f$ at~$x$ assumes the existence of every lower order Peano derivative of~$f$ at $x$.
Moreover, it is well known that, by Taylor expansion about $x$, the existence of the~$n$th Peano
derivative $f_{\left(  n\right)  }(x)$ forces every generalized Riemann
derivative of order $n$ to exist and to be equal to $f_{\left(  n\right)
}\left(  x\right)  $.

Riemann derivatives were introduced by Riemann in the mid 1800s; see \cite{R}. Generalized Riemann derivatives were introduced by Denjoy in \cite{D} in 1935; they have many applications in the theory of trigonometric series~\cite{SZ,Z}. Generalized Riemann derivatives have been shown to satisfy basic properties of ordinary derivatives, such as monotonicity, convexity, or the mean value theorem. \cite{AJ,FFR,GGR1,HL,HL1,MM,T,W} Multidimensional Riemann derivatives were studied in \cite{AC1}, and quantum Riemann derivatives appeared in~\cite{AC,ACR}.
Peano derivatives are due to Peano in \cite{P} in 1892 and were further developed by de la Vall\'ee Poussin in \cite{VP} in 1908. For more recent developments on Peano derivatives, see~\cite{F,F1,FR,FW,LPW}.
Surveys on generalized Riemann and Peano derivatives are found in \cite{As2,EW,O}.

\subsection{Context and Motivation}The main problem in the area of generalized derivatives has been the equivalence between Peano and generalized Riemann derivatives, that is, the problem of which $\mathcal{A}$ has the property that the existence of the generalized Riemann derivative $D_{\mathcal{A}}f(x)$ is equivalent to the existence of the~$n$th Peano derivative $f_{(n)}(x)$, for all functions~$f$ at~$x$. This problem goes back to 1927, when Khintchine proved in \cite{Ki} that the first symmetric derivative~$D_1^sf(x)$ is equivalent to the first Peano derivative $f_{(1)}(x)$ at~a.e. points on a measurable set. This was extended greatly by Marcinkiewicz and Zygmund in \cite{MZ} to the~$n$th symmetric Riemann derivative $D_n^sf(x)$ being a.e. equivalent to~$f_{(n)}(x)$, and further generalized by Ash in \cite{As} to each $n$th generalized Riemann derivative $D_{\mathcal{A}}f(x)$ is a.e. equivalent to the $n$th Peano derivative $f_{(n)}(x)$ on a measurable set.~The measurability condition was removed by Fejzi\'c and Weil in \cite{FW}.
The problem was completely solved in 2017 by Ash, Catoiu, and Cs\"ornyei with the result of \cite[Theorem 1]{ACCs}, asserting that in orders $n\geq 2$ no generalized Riemann derivative is equivalent to the $n$th Peano derivative for all functions at $x$, and in order $n=1$ a non-trivial complete set of generalized Riemann derivatives equivalent to the first Peano derivative is provided.

In this paper, we look at the above largely negative resolution to the main problem and suggest that, in orders $n\geq 2$, the right problem to look at, that has a positive outcome, is not the equivalence between Peano and a single generalized Riemann derivative, but the equivalence between Peano and sets of generalized Riemann derivatives, a new direction of research in generalized derivatives. 

The motivation for this comes from an earlier result and conjecture by Ginchev, Guerragio and Rocca in \cite{GGR}, see also \cite{GR}, from 1998, predicting the existence of such a set of generalized Riemann derivatives. Their conjecture is that the existence of  $D_1f(x)$ and each of the~$n(n-1)/2$ derivatives,\vspace{-.1in}
\[
D_{k,-j}f(x)=\lim_{h\rightarrow 0} h^{-k}\sum_{i=0}^k(-1)^i\binom ki f(x+(k-j-i)h),\vspace{-.05in}
\]
for $k=2,\ldots,n$ and $j=0,\ldots ,k-2$, is equivalent to the existence of the $n$th Peano derivative $f_{(n)}(x)$, for all $f$ at $x$. They proved the predicted result by hand for $n=1,2,3,4$ and with computer assistance for $n\leq 8$, leaving the cases $n>8$ as a conjecture.

The goal of this paper is not to prove the GGR Conjecture, but a variant of it that has all of $D_{k,-j}$ replaced by $D_{k,j}$, for $j$ in the same range. This is done in Theorem~\ref{T2}(iv), providing a set of generalized Riemann derivatives equivalent to the~$n$th Peano derivative for all functions at $x$. In proving this theorem, we employ two technical results: an old lemma due Marcinkiewicz and Zygmund in \cite{MZ}, which we list here as Theorem~\ref{T2}(i) and refer to it as Theorem~MZ, that has been known to be useful in similar situations; and a new result of our own, Corollary~\ref{C1} below. The connection between these two ingredients led to the result of Theorem~\ref{T2}(ii), providing another set of generalized Riemann derivatives equivalent to the $n$th Peano derivative for all $f$ at $x$. These are our two pointwise characterizations of the Peano derivative by sets of generalized Riemann derivatives, and the beginning of a new subject. More details about this subject are given in the note at the end of the article.

\subsection{Results}
Even before we have a set of generalized Riemann derivatives equivalent to the $n$th Peano derivative for all functions at a point, the following lemma and its corollary provide a dual way of looking at such sets. The lemma has a straightforward inductive proof, and the corollary is an immediate consequence of the lemma. We skip both of their proofs.
\begin{lem}
Given a positive integer $n$, for each $k$ with $1\leq k\leq n$, let $\mathscr{D}_k$ be a finite nonempty set of $k$th generalized Riemann derivatives such that $\mathscr{D}_1=\{D_1\}$. Then the following statements are equivalent for all functions $f$ at $x$:
\begin{enumerate}
\item[(i)\,] $(\mathscr{D}_1\cup \mathscr{D}_n\cup\ldots\cup\mathscr{D}_n)f(x)$ exists $\Longleftrightarrow $ $f_{(n)}(x)$ exists.
\item[(ii)] For $k=2,\ldots,n$, both $f_{(k-1)}(x)$ and $\mathscr{D}_kf(x)$ exists $\Longleftrightarrow $ $f_{(k)}(x)$ exists.
\end{enumerate}
\end{lem}

\begin{corollary}\label{C1}
For each positive integer $n$, let $\mathscr{D}_n$ be a finite nonempty set of $n$th generalized Riemann derivatives such that $\mathscr{D}_1=\{D_1\}$. Then the following statements are equivalent for all functions $f$ at $x$:
\begin{enumerate}
\item[(i)\,] For all $n\geq 2$, $(\mathscr{D}_1\cup \mathscr{D}_2\cup\ldots\cup\mathscr{D}_n)f(x)$ exists $\Longleftrightarrow $ $f_{(n)}(x)$ exists.
\item[(ii)] For all $n\geq 2$, both $f_{(n-1)}(x)$ and $\mathscr{D}_nf(x)$ exists $\Longleftrightarrow $ $f_{(n)}(x)$ exists.
\end{enumerate}
\end{corollary}
Our main result is the following theorem, providing the two announced characterizations of the $n$th Peano derivative by sets of generalized Riemann derivatives.
Let $\widetilde{D}_nf(x)$ be the exact generalized Riemann derivative based at $x,x+h,x+2h,x+2^2h,\ldots ,x+2^{n-1}h$, and let ${D}_n^{\rm sh}f(x)=\{D_{n,j}f(x)\,|\,j=0,1,\ldots , n-2\}$ be the set of the first $n-1$ forward shifts of the $n$th forward Riemann derivative $D_nf(x)$.

\begin{theorem}[Two pointwise characterizations of the Peano derivative]\label{T2}
If $n\geq 2$, then,
for all functions $f$ and points $x$,

\begin{tabular}{ll}
{\rm (i)} {\rm (}{Theorem MZ}{\rm )} both $f_{(n-1)}(x)$ and $\widetilde{D}_nf(x)$ exist &$\Longleftrightarrow $ $f_{(n)}(x)$ exists;\\
{\rm (ii)} all $\widetilde{D}_1f(x),\widetilde{D}_2f(x),\ldots ,\widetilde{D}_nf(x)$ exist & $\Longleftrightarrow $ $f_{(n)}(x)$ exists;\\
{\rm (iii)} both $f_{(n-1)}(x)$ and ${D}_n^{\rm sh}f(x)$ exist &$\Longleftrightarrow $ $f_{(n)}(x)$ exists;\\
{\rm (iv)} $D_1f(x)$ and all ${D}_2^{\rm sh}f(x),\ldots ,{D}_n^{\rm sh}f(x)$ exist & $\Longleftrightarrow $ $f_{(n)}(x)$ exists.
\end{tabular}
\end{theorem}

The reverse implications in all parts of the above theorem come as consequences of the property of the $n$th Peano derivative, that it implies every lower order Peano derivative and every $n$th generalized Riemann derivative, for all functions at $x$.

Part (i) of the above theorem is the result of Theorem~MZ, due to Marcinkiewicz and~Zygmund in \cite[Lemma~1]{MZ}.
Part (ii) is our first pointwise characterization of the $n$th Peano derivative by sets of generalized Riemann derivatives, for all functions~$f$ at $x$. By Corollary~\ref{C1}, parts (i) and (ii) are equivalent to each other, and so, part~(i) also represents the first characterization of the $n$th Peano derivative by sets of generalized Riemann derivatives, this time for all~$n-1$ times Peano differentiable functions~$f$ at~$x$. Indeed, Theorem~MZ can be stated as the $n$th generalized Riemann derivative~$\widetilde{D}_nf(x)$ is equivalent to the $n$th Peano derivative for all~$n-1$ times Peano differentiable functions~$f$ at~$x$, an interpretation that Marcinkiewicz and Zygmund did not envision for their lemma.

Similarly, parts (iii) and (iv) are also equivalent to each other; they represent our second pointwise characterization of the Peano derivative. The difference between them is that part (iv) is the second characterization of the $n$th Peano derivative by sets of generalized Riemann derivatives, for all functions at~$x$; while part (iii) is the second characterization of the $n$th Peano derivative by sets of generalized Riemann derivatives, for all $n-1$ times Peano differentiable functions at $x$.

Based on the known proof of part (i) and the equivalences between the two pairs of parts, the proof of Theorem~\ref{T2} is reduced to the proof of the direct implication part~(iii).~The entire paper is devoted to this proof.
In Section~\ref{S1}, the proof is further reduced to a result on recursive sets, Lemma~\ref{P4}, which is then proved in Section~\ref{S2} using an interesting combinatorial algorithm.

The use of a combinatorial method of proof here should not come as a surprise, since the expression of each forward shift of a Riemann derivative involves binomial coefficients, and, as it was recently proved in \cite{ACF1}, the expressions of the special generalized Riemann derivatives $\widetilde{D}_nf(x)$ involve $q$-binomial coefficients, for $q=2$.
\[
\ast\quad\ast\quad\ast
\]
Before proceeding with the proof of the remaining part of the main theorem, the direct implication in Theorem~\ref{T2}(iii), our last result is meant to help with understanding the limits of Theorem~MZ. It shows that, for~$n\geq 3$, the two most widely known $n$th generalized Riemann derivatives, the $n$th symmetric Riemann derivative $D_n^sf(x)$ and the $n$th Riemann derivative~$D_nf(x)$, fail to characterize the $n$th Peano derivative $f_{(n)}(x)$ at $x=0$ in the same way as $\widetilde{D}_nf(x)$ does in Theorem MZ.

\begin{theorem}\label{T0} Suppose $n\geq 3$. Then there exist functions $f$ and real numbers $x$ such that, separately:

{\rm (i)} both $f_{(n-1)}(x)$ and $D_n^sf(x)$ exist, but $f_{(n)}(x)$ does not exist;

{\rm (ii)} both $f_{(2)}(x)$ and $D_3f(x)$ exist, but $f_{(3)}(x)$ does not exist.
\end{theorem}

\begin{proof}
(i) Consider the functions $g:[0,\infty )\rightarrow \mathbb{R}$, defined by $g(x)=x^{n-1/2}$, and $f:\mathbb{R}\rightarrow \mathbb{R}$, defined by
\[
f(x)=\begin{cases}
g(x) &\text{, if }x\geq 0,\\
(-1)^{n-1}g(-x)&\text{, if }x<0.
\end{cases}
\]
Clearly, $g(h)=o(h^{n-1})$ implies $f(h)=o(h^{n-1})$, so $f$ is $n-1$ times Peano differentiable at $0$ and $f_{(0)}(0)=f_{(1)}(0)=\cdots =f_{(n-1)}(0)=0$, while $\lim_{h\rightarrow 0^+} g(h)/h^n=\lim_{h\rightarrow 0^+}1/\sqrt{h}=\infty $ implies that $f_{(n)}(0)$ does not exist. On the other hand, since~$f$ is of opposite parity as $n$ is, the $n$th symmetric derivative $D_n^sf(0)$ exists and is equal to zero.

(ii) Let $G=\langle 2,3\rangle $ be the multiplicative subgroup of the rationals generated by the integers 2 and 3. Then $G=\{ 2^m3^n\mid m,n\text{ integers}\}$, so that if $h\in G$, then $2h,3h\in G$ and if $h\notin G$ then $2h,3h\notin G$. Define $f\colon \mathbb{R}\rightarrow \mathbb{R}$ by 
\[
f(x)=\begin{cases}
(-1)^{m+n}x^3 &\text{, if }x={2^m}{3^n}\in G,\\
0&\text{, if }x\notin G.
\end{cases}
\]

Since $0\leq |f(x)|\leq x^3$, $f$ is continuous at 0, has two Peano derivatives at 0, and $f_{(0)}(0)=f_{(1)}(0)=f_{(2)}(0)=0$. The third Peano derivative $f_{(3)}(0)$ does not exist, since the ratio $f(h)/h^3$ has three distinct limit points, 0, 1 and -1, as $h\rightarrow 0$.

Let $h=2^m3^n$. Then $\Delta_3(0,h;f)=f(3h)-3f(2h)+3f(h)-f(0)$
\[
\begin{aligned}
&=(-1)^{m+n+1}3^3h^3-3(-1)^{m+n+1}2^3h^3+3(-1)^{m+n}h^3-0\\
&=(-1)^{m+n}h^3(-27+24+3)\\
&=0,
\end{aligned}
\]
so $f$ is three times forward Riemann differentiable at 0 and $D_3f(0)=0$.
\end{proof}

We suspect the result in part (ii) of Theorem~\ref{T0} to hold for a general $n$ in place of $n=3$. Since the example we use here does not seem to extend to higher $n$, and such an extension does not seem to be of as much interest as the remainder of this paper, we leave this extension as an open question.
However, Theorem~\ref{T0} as stated fulfills its principal goal of motivating the main theorem of the paper, Theorem~\ref{T2}, as much as its most general version would.
This conjecture is the natural extension to sets of generalized Riemann derivatives of the much simpler statement that, in orders $n\geq 2$, the Riemann derivative $D_nf(x)$ is not equivalent to the Peano derivative $f_{(n)}(x)$ for all functions $f$ at $x$. Had this been stated by Peano in 1892, it would have been the longest lasting solved open problem in generalized differentiation. Its solution came recently as a part of the main result in \cite{ACCs}.

For more on the recent developments on the subject we started in here, as well as on the importance of the article, see the note at the end.

\section{Main theorem reduced to a combinatorial problem}\label{S1}

To prove the direct implication in Theorem~\ref{T2}(iii), by the direct implication in Theorem~MZ, it suffices to show that

\begin{lemma}\label{L3}
For each function $f$ and real number $x$,
\[
{D}_n^{\rm sh}f(x)\text{ exists }\Longrightarrow \widetilde{D}_nf(x)\text{ exists.}
\]
\end{lemma}

Note that, unlike Theorem~\ref{T2}(iii), the above lemma relates only $n$th generalized Riemann derivatives, without involving any Peano derivative.

Before proceeding with the proof of Lemma~\ref{L3}, it is important to lay down a few ideas about generating new generalized Riemann derivatives from old ones, which we will refer to as \emph{operations} with generalized Riemann derivatives. Specifically, we will concentrate on two basic operations with $n$th generalized Riemann derivatives based at $n+1$ points:

\subsubsection*{1. Dilation} The dilation by a real number $r$ of an $n$th generalized Riemann difference
$\Delta_{\mathcal{A}}(x,h;f)=\sum_{i=0}^nA_if(x+a_ih)$ corresponding to $\mathcal{A}=\{A_i;a_i\,|\,i=0,\ldots,n\}$ is the~$n$th generalized Riemann difference,
\[\Delta_{\mathcal{A}_r}(x,h;f)=\sum_{i=0}^nr^{-n}A_if(x+a_irh),\]
corresponding to the vector $\mathcal{A}_r=\{r^{-n}A_i;ra_i\,|\,i=0,\ldots,n\}$. Moreover, a function $f$ is $\mathcal{A}$-differentiable at $x$ if and only if is $\mathcal{A}_r$-differentiable at $x$ and $D_{\mathcal{A}}f(x)=D_{\mathcal{A}_r}f(x)$.

To simplify the language, for fixed $f$ and $x$, we say that an ordered set $\mathcal{A}$ of $2n+2$ elements satisfying the $n$th Vandermonde conditions is \emph{good} if $f$ is $\mathcal{A}$-differentiable at~$x$. Then $\mathcal{A}$ is good if and only if $\mathcal{A}_r$ is good. And since each $n$th generalized Riemann derivative based at $n+1$ points is uniquely determined by its base points, $\mathcal{A}$ is uniquely determined by $\{a_0,\ldots ,a_n\}$, while $A_0,\ldots ,A_n$ are just place holders. In this way, a plain set $\{a_0,\ldots ,a_n\}$ is good if and only~if, for each non-zero real number~$r$, its \emph{$r$-dilate} $\{ra_0,\ldots ,ra_n\}$ is good. 

\subsubsection*{2. Elimination} By the $n$th Vandermonde conditions, only a unique non-zero scalar multiple of each non-zero linear combination of two $n$th generalized Riemann differences is a generalized Riemann difference. The set of base points of such a non-zero linear combination is the union of the base points of the terms, minus the base points corresponding to the terms that got eliminated by the linear combination. Focusing only on pairs of $n$th generalized Riemann differences based at $n+1$ points and having the additional property that they share $n$ base points, observe that by taking non-zero linear combinations of the differences corresponding to such pairs, one arrives at one of the following three situations:
\begin{itemize}
\item No common terms of the two differences were eliminated by the linear combination. Then the resulting $n$th difference has $n+2$ base points, a discarded case.
\item At least two common terms were eliminated by the linear combination. Then the resulting $n$th difference has $\leq n$ base points, an impossibility.
\item Only one common term of the two differences got eliminated by the linear combination. Then the resulting $n$th difference has $n+1$ base points, hence it is of the desired kind. And since linear combinations of $n$th generalized differences are scalar multiples of $n$th generalized differences, if $S=\{a_0,\ldots ,a_n\}$ and $T=\{b_0,\ldots ,b_n\}$ are good sets such that $|S\cap T|=n$ then, for each $a\in S\cap T$, the set $S\cup T\setminus \{a\}$ is also a good set.
\end{itemize}
As an example highlighting the third bullet situation, consider the two shifts,
\[
\begin{aligned}
\Delta_{3,0}(h)&=f(x+3h)-3f(x+2h)+3f(x+h)-f(x)\text{ and }\\
\Delta_{3,1}(h)&=f(x+4h)-3f(x+3h)+3f(x+2h)-f(x+h),
\end{aligned}
\]
of the third Riemann difference $\Delta_3(h)$ of $f$ at $x$, so that the sets $\{0,1,2,3\}$ and $\{1,2,3,4\}$ representing their base points are good.
The linear combination $\tfrac 34\Delta_{3,0}(h)+\tfrac 14\Delta_{3,1}(h)$
that eliminates the term in $f(x+3h)$ is the third exact generalized Riemann difference corresponding to the (new good) set $\{0,1,2,4\}$ obtained by eliminating the entry 3 between the other two known good sets.

Back to Lemma~\ref{L3}, its statement in language of good sets goes as follows: If the sets $\{0,1,2,\ldots ,n\},\{1,2,3,\ldots ,n+1\},\ldots ,\{n-2,n-1,n,\ldots ,2n-2\}$ are good, then so is the set $\{0,1,2,4,\ldots ,2^{n-1}\}$.

As for its proof, since we already know two operations that produce new good sets from old, dilation and elimination, to prove the lemma, it suffices to provide an algorithm that inputs the given sets in the above hypothesis and outputs the set in the conclusion, by only using dilations of given or previously deduced sets, and elimination of a common element between a pair of given or deduced sets that have $n$ common elements.

\medskip
Summarizing, we have reduced the proof of Lemma~\ref{L3}, and implicitly the one of Theorem~\ref{T2}(iii), to the following result of recursive set theory.

\begin{lemma}\label{P4}
Let $n$ be an integer at least 2. A collection $\mathcal{S}$ of sets, each consisting of $n+1$ non-negative integers, is defined by the following properties:
\begin{enumerate}
\item[({\rm i})] $\{0,1,2,\ldots ,n\},\{1,2,3,\ldots ,n+1\},\ldots ,\{n-2,n-1,n,\ldots ,2n-2\}\in\mathcal{S}$;
\item[({\rm ii})] if $S \in \mathcal{S}$, then $2S:=\{2s\,|\,s\in S\} \in \mathcal{S}$;
\item[({\rm iii})] if $S,T \in \mathcal{S}$ have $|S\cap T|=n$, then for each $a\in S\cap T$, $S\cup T\setminus \{a\}\in \mathcal{S}$.
\end{enumerate}
Then $\{0,1,2,4,\ldots ,2^{n-1}\}\in \mathcal{S}$.
\end{lemma}

\section{Proof of the combinatorial problem}\label{S2}

The combinatorial Lemma~\ref{P4}, whose complex proof lies at the heart of this paper, relies on a four-step combinatorial algorithm on recursive sets, described in Section 2.2. All ideas behind the algorithm are brought up earlier in Section~2.1, where several smaller cases are investigated.

As a simplifying terminology for the proof,
we convene, in all displayed equations, to write each ordered set as a row-vector, that is, without braces and commas.
The hypothesis~(i) means \emph{input} the given sets (vectors) of $\mathcal{S}$ into the algorithm, (ii) is the \emph{dilation} of a set by 2, and (iii) is the \emph{elimination} of a common element between two sets that share $n$ elements. And two elements $S,T$ of  $\mathcal{S}$, with $|S\cap T|=n$, are said to be \emph{set for elimination}.

\subsection{Smaller cases}\label{S2.1}
In this subsection we prove Lemma~\ref{P4}, for $n=3,4,7,10$. Each of these cases contains an extra idea for the general proof that was not shown in the previous cases, so that the case $n=10$ has all ideas needed for the general proof.

The proof of the $n=3$ case of Lemma~\ref{P4} was accomplished at the end of the previous section by displaying $\frac 34\Delta_{3,0}(h)+\frac 14\Delta_{3,1}(h)$. As a first example, we rewrite the $n=3$ case here to begin illustrating the notation and language we'll use in giving the general proof of Lemma~\ref{P4}:
Input the two given elements $\{0,1,2,3\},\{1,2,3,4\}\in \mathcal{S}$ as row vectors, with their common elements listed under each other, and highlight the number 3 as the only non-zero entry that is not a power of 2, which we will refer to as an \emph{intruder}.
\[
\begin{array}{ccccc}
0&1&2&\mathbf{3}&\\&1&2&\mathbf{3}&4
\end{array}
\]
Eliminate 3 between the two given sets to deduce that $\{0,1,2,4\}\in \mathcal{S}$, as needed.

\subsection*{Case $n=4$} Write the three given sets as row-vectors in a parallelogram array, with the common elements written in the same column, and mark all intruders.
\[
\begin{array}{ccccccc}
0&1&2&\mathbf{3}&4&&\\&1&2&\mathbf{3}&4&\mathbf{5}&\\&&2&\mathbf{3}&4&\mathbf{5}&\mathbf{6}
\end{array}
\]
Note that consecutive rows are set for elimination. This property will stay in place on all arrays till the end of the algorithm.

Consider the following chain of eliminations.
Replace the third row with the result of eliminating 5 between itself and the second row, to deduce that $\{1,2,3,4,6\}\in\mathcal{S}$.   Then replace this new third row with the result of eliminating 3 between itself and the first row, to deduce that $\{0,1,2,4,6\}\in \mathcal{S}$. The conjunction of these two eliminations is read as follows: on the row ending in 6 (even number), we eliminated as many odd numbers (two, namely 5 and 3, in this order) as there are rows above it, at the price of adding 0 and 1 at the beginning of the row.

Our first processing of the above parallelogram array is the set of all possible elimination compounds described above: Starting at the bottom and going up, on each row ending with an even entry (first and third), eliminate as many odd entries (high to low) as there are rows above it.
Then delete the remaining (second) row(s), as shown on the left chart below. Then double the first row and move it to the bottom,
\[
\begin{array}{cccccc}
0&1&2&\mathbf{3}&4&\\0&1&2& &4& \mathbf{6}
\end{array}
\qquad \qquad \qquad
\begin{array}{cccccc}
0&1&2 &4 &\mathbf{6}&\\0&&2&4&\mathbf{6}&8
\end{array}
\]
as shown on the right chart. Finally, eliminate the last intruder $6$ between the two rows to deduce that $\{0,1,2,4,8\}\in \mathcal{S}$, as needed.

As a shortcut for the cojunction of the last two operations, we say that the same conclusion was obtained directly from the earlier array by eliminating 6 between the second row and twice the first row. In this way, for a general $n$, from this stage on to the end of the algorithm, all arrays will have the additional property that twice the top row and the bottom row are set for elimination.

\subsection*{Case $n=7$} Input the given vectors in a, by now familiar, parallelogram array.
\[
\begin{array}{ccccccccccccccc}
0&1&2&\mathbf{3}&4&\mathbf{5}&\mathbf{6}&\mathbf{7}&&&&&\\
&1&2&\mathbf{3}&4&\mathbf{5}&\mathbf{6}&\mathbf{7}&8&&&&\\
&&2&\mathbf{3}&4&\mathbf{5}&\mathbf{6}&\mathbf{7}&8&\mathbf{9}&&&\\
&&&\mathbf{3}&4&\mathbf{5}&\mathbf{6}&\mathbf{7}&8&\mathbf{9}&\mathbf{10}&&\\
&&&&4&\mathbf{5}&\mathbf{6}&\mathbf{7}&8&\mathbf{9}&\mathbf{10}&\mathbf{11}&\\
&&&&&\mathbf{5}&\mathbf{6}&\mathbf{7}&8&\mathbf{9}&\mathbf{10}&\mathbf{11}&\mathbf{12}
\end{array}
\]
Eliminate as many odd entries, high to low, from the rows ending in even entries as there are rows above them, and delete all rows ending in odd entries.
\[
\begin{array}{cccccccccccc}
0&1&2&\mathbf{3}&4&\mathbf{5}&\mathbf{6}&8&&\\
0&1&2&\mathbf{3}&4&&\mathbf{6}&8&\mathbf{10}&\\
0&1&2&&4&&\mathbf{6}&8&\mathbf{10}&\mathbf{12}
\end{array}
\]
Note that now all rows have the same number of intruders: three. Replace the third row with the result of the elimination of $12$ between itself and the double of the first row. Replace the second row with the result of the elimination of $10$ between itself and the newly established third row. Replace the first row with the result of the elimination of $6$ between itself and the new second row. The effect of this process is to replace the largest intruder in each row with 16, the next power of $2$. We refer to this chain of actions as~\emph{cutting off the largest intruders} in each row.
\[
\begin{array}{cccccccccccc}
0&1&2&\mathbf{3}&4&\mathbf{5}&&8&&16\\
0&1&2&\mathbf{3}&4&&\mathbf{6}&8&&16\\
0&1&2&&4&&\mathbf{6}&8&\mathbf{10}&16
\end{array}
\]
Repeat cutting off intruders from all rows, except for the first one, which is deleted. The difference between this cut off that deletes the first row and the preceding one that retained its first row is marked by the last intruder on the first row now being odd 5, instead of last time being even 6.
\[
\begin{array}{cccccccccc}
0&1&2&\mathbf{3}&4&&8&16&32\\
0&1&2&&4&\mathbf{6}&8&16&32
\end{array}
\]
Finally, since the last intruder in the first row is the odd number 3, we delete the first row and cut off the last intruder in the second row to get $\{0,1,2,4,8,16,32,64\}\in\mathcal{S}$.

\subsection*{Case $n=10$} Start with the parallelogram made with the row-vectors corresponding to the given sets
$\{0,1,\ldots,10\},\{1,2,\ldots ,11\},\ldots ,\{8,9,\ldots,18\}$ of $\mathcal{S}$, where each column has equal entries, and from the bottom to the top, in each row ending with an even entry, eliminate as many of the largest odd intruders as there are rows above it. Moreover, eliminate the rows ending in odd entries.
\[
\begin{array}{cccccccccccccll}
0&1&2&\mathbf{3}&4&\mathbf{5}&\mathbf{6}&\mathbf{7}&8&\mathbf{9}&\mathbf{10}&&&&\\
0&1&2&\mathbf{3}&4&\mathbf{5}&\mathbf{6}&\mathbf{7}&8&&\mathbf{10}&\mathbf{12}\\
0&1&2&\mathbf{3}&4&\mathbf{5}&\mathbf{6}&&8&&\mathbf{10}&\mathbf{12}&\mathbf{14}\\
0&1&2&\mathbf{3}&4&&\mathbf{6}&&8&&\mathbf{10}&\mathbf{12}&\mathbf{14}&16^{\ast }\\
0&1&2&&4&&\mathbf{6}&&8&&\mathbf{10}&\mathbf{12}&\mathbf{14}&16&\mathbf{18}
\end{array}
\]
The largest numbers in each row are the even numbers between $n=10$ and $2n-2=18$. There is a single power of $2$ among them, namely $16$, which we mark by an asterisk. Note that the number of intruders in each of the rows above the asterisk is six, while the rows with asterisk or below it have only five intruders. Using the asterisk row as a base, cut off the largest intruders in all rows above it.
\[
\begin{array}{ccccccccccccclc}
0&1&2&\mathbf{3}&4&\mathbf{5}&\mathbf{6}&\mathbf{7}&8&\mathbf{9}&&&&16&\\
0&1&2&\mathbf{3}&4&\mathbf{5}&\mathbf{6}&\mathbf{7}&8&&\mathbf{10}&&&16\\
0&1&2&\mathbf{3}&4&\mathbf{5}&\mathbf{6}&&8&&\mathbf{10}&\mathbf{12}&&16\\
0&1&2&\mathbf{3}&4&&\mathbf{6}&&8&&\mathbf{10}&\mathbf{12}&\mathbf{14}&16^{\ast }\\
0&1&2&&4&&\mathbf{6}&&8&&\mathbf{10}&\mathbf{12}&\mathbf{14}&16&\mathbf{18}
\end{array}
\]
Now there are five intruders in each row. As long as the last entry in the first row is odd, delete it and (using its double, which we do not write down, as a base) eliminate the highest intruders in the remaining rows. We deduce

\[
\begin{array}{cccccccccccclc}
0&1&2&\mathbf{3}&4&\mathbf{5}&\mathbf{6}&\mathbf{7}&8&&&&16&32\\
0&1&2&\mathbf{3}&4&\mathbf{5}&\mathbf{6}&&8&\mathbf{10}&&&16^{\ast }&32\\
0&1&2&\mathbf{3}&4&&\mathbf{6}&&8&\mathbf{10}&\mathbf{12}&&16&32\\
0&1&2&&4&&\mathbf{6}&&8&\mathbf{10}&\mathbf{12}&\mathbf{14}&16&32
\end{array}
\]
and
\[
\begin{array}{cccccccccclll}
0&1&2&\mathbf{3}&4&\mathbf{5}&\mathbf{6}&8&&&16^{\ast }&32&64\\
0&1&2&\mathbf{3}&4&&\mathbf{6}&8&\mathbf{10}&&16&32&64\\
0&1&2&&4&&\mathbf{6}&8&\mathbf{10}&\mathbf{12}&16&32&64
\end{array}
\]
When the last intruder in the first row is even, we keep the first row and cut off the largest intruders from all rows. 
\[
\begin{array}{ccccccccllcccc}
0&1&2&\mathbf{3}&4&\mathbf{5}&&8&&16^{\ast }&32&64&128\\
0&1&2&\mathbf{3}&4&&\mathbf{6}&8&&16&32&64&128\\
0&1&2&&4&&\mathbf{6}&8&\mathbf{10}&16&32&64&128
\end{array}
\]
Repeat the process until all intruders are eliminated. We get
\[
\begin{array}{cccccccccccccc}
0&1&2&\mathbf{3}&4&&8&16&32&64&128&256\\
0&1&2&&4&\mathbf{6}&8&16&32&64&128&256
\end{array}
\]
and, finally, deduce that $\{0,1,2,4,8,16,32,64,128,256,512\}\in \mathcal{S}$, as desired.

\subsection{Proof of Lemma~\ref{P4}}\label{S2.2} The proof of the general case in Lemma~\ref{P4} follows ideas from the proofs of the smaller cases described in the first half of the section. The general proof is based on the following combinatorial (Gaussian) elimination algorithm:

\subsubsection*{STEP 1} Arrange (input) the given sets
\[\{0,1,\ldots ,n\},\{1,2,\ldots ,n+1\},\ldots ,\{n-2,n-1,\ldots,2n-2\}\]
in order, one under the other, as row vectors in a parallelogram array so that each next vector is shifted one position to the right of the previous vector, so that all equal entries are a part of the same column. Highlight all intruders.

\subsubsection*{STEP 2} Replace each row ending in an even entry with the one obtained from it by deleting as many of its largest odd entries as there are rows above it and adding all non-negative integers smaller than its smallest entry. Delete all rows ending in an odd entry.

\medskip
In this way, the new $k$th row from the bottom, for $k=1,2,\ldots,\lfloor (n-1)/2\rfloor $, is coming from the old row
$\{n-2k,n-2k+1,n-2k+2,\ldots ,2(n-k)\}$
by inserting all non-negative integers, $0,1,\ldots,n-2k-1$, that can fit in front of it, and deleting its largest $n-2k$ odd entries $2(n-k)-1,2(n-k)-3,\ldots ,2(n-k)-2(n-2k)+1=2k+1$. The new $k$th row (set) from the bottom has the expression:
\begin{equation}\label{st2}
\{0,1,\ldots ,2k-1,2k,2(k+1),2(k+2),\ldots ,2(n-k-1),2(n-k)\}.
\end{equation}
Then the top row is either $\{0,1,\ldots,n\}$, for even $n$, or $\{0,1,\ldots,n-2,n-1,n+1\}$, for odd $n$, and the bottom row is $\{0,1,2,4,6,8,\ldots ,2(n-1)\}$. Note that all consecutive rows are set for elimination. This will continue to the end of the algorithm.

\subsubsection*{STEP 3} Mark with an asterisk the unique power of 2 among the last entries in all rows. Note that every row strictly above this unique row has the same number of intruders, and the rest of the rows all have one fewer than this common number. Cut off the largest intruders in the rows above the asterisk. The effect is that these will be replaced by the power of 2 marked with asterisk, and all rows will then have the same number of intruders.

\medskip
Suppose the asterisk entry $2^{\eta }$ occurs in the $\ell $th row. Then after Step 3,
the $k$th row from the bottom, for $k\leq \ell $, is as in (\ref{st2}), while the one for $k>\ell $ is
\begin{equation*}\label{st3}
\{0,1,\ldots ,2k-1,2k,2(k+1),2(k+2),\ldots ,2(n-k-1),2^{\eta }\}.
\end{equation*}
The top row is $\{0,1,\ldots,n-1,2^{\eta }\}$, regardless of the parity of $n$, and the bottom row is $\{0,1,2,4,6,8,\ldots ,2(n-1)\}$.
So twice the top row and the bottom row are set for elimination, and the next step is granted.

\subsubsection*{STEP 4} If the last intruder in the top row is even, cut off the largest intruders in all rows. Otherwise, delete the first row and cut off  the largest intruders in the remaining rows.  

\smallskip
Repeat Step 4 until all intruders are eliminated and the array has shrunk to a single row,
the desired $\{0,1,2,4,8,\ldots ,2^{n-1}\}$.

This completes the proof of Lemma~\ref{P4}. Furthermore, this completes the proof of Theorem~\ref{T2}, thereby achieving the main goal of the paper.

\subsection*{Remarks} We end by giving a few more details to assist the reader's comprehension of how each application of either of the last two steps in the algorithm is granted by the conclusion of the previous step.

After Step 2, the odd intruders in the $k$th row from the bottom are $3,5,\ldots ,2k-1$. Their number is $\alpha (k)=k-1$. The even intruders in the same row are all positive even entries, minus all powers of 2. Their count is $\beta (k)=n-k-\eta (k)$, where $\eta (k)=\lfloor \log_22(n-k)\rfloor $. Then $\eta (k)=\eta $, when $k\leq \ell $, and $\eta (k)=\eta +1$, when $k>\ell $. The total number of rows is $\nu =\lfloor n/2\rfloor $, and the total number of intruders in row $k$, denoted as $\gamma (k)$, is $n-\eta $, when $k>\ell $, and $n-\eta -1$, when $k\leq \ell $. Thus all actions in Step 3 are granted.

At the end of Step 3, or before each application of Step 4,
if the last intruder in the top row is even, then Step 4 removes it together with its double in the bottom row, so that the new top and bottom rows are set for elimination, based on the same property for the old rows. If the last intruder in the first row is odd, say $2s+1$, then the top row is $0,1,\ldots ,2s+1$ followed by 2-powers, the row below it is $0,1,\ldots ,2s$, and one more even intruder followed by 2-powers, and the bottom row is $0,1,2,4,6,\ldots , 4s,4s+2$ followed by 2-powers. After Step 4, the top row will be $0,1,\ldots ,2s$ followed by 2-powers, and the bottom row will be $0,1,2,4,6,\ldots , 4s$ followed by 2-powers, making them set for elimination. Thus the next application of Step 4 is granted in all cases.

\subsection*{Note.} This is the first article attacking the general proof of the 1998 conjecture by Ginchev, Guerragio, and Rocca, by proving a variant of the conjecture in Theorem~\ref{T2}(iv). Unfortunately, due to delays in the publication process, a number of follow-up papers have already appeared in print since the article was completed: (1) The GGR Conjecture has been solved in~\cite{AC2}, with a proof that relies in part on the result of Lemma~\ref{P4} and its proof; in particular, in addition to its own merits, the present article is necessary to validate this proof. (2)~A~new proof of Lemma~\ref{L3} is given in~\cite{ACF1}. (3)~A~result directed towards extending Theorem~\ref{T2}(iii) for $n=1$ appeared in \cite{C}. (4) A generalization of the GGR Conjecture, for generalized Riemann smoothness instead of generalized Riemann differentiation, is proved in \cite{CF}. (5) More sets of generalized Riemann derivatives proven to be equivalent to the $n$th Peano derivative for all functions at $x$ are found in \cite{ACF}. (6) The same article also confirmed the $n=7$ case of our textual conjecture given at the end of the introduction, that the result of Theorem~\ref{T0}(ii) extends from $n=3$ to general $n$. (7) There has been an effort toward extending the result of Theorem~\ref{T2} to more general sets of generalized Riemann derivatives. In particular, it was conjectured in \cite{ACF} that the Gaussian derivatives highlighted in there are all such derivatives that can play the role of $\widetilde{D}_nf(x)$ in Theorem~MZ. A number of counterexamples given in \cite{ACF2} have shown that the problem of finding all generalized Riemann derivatives that are equivalent to the $n$th Peano derivative for all $n-1$ times Peano differentiable functions at a point is a much harder problem. Consequently, that article elevated our textual conjecture as the most outstanding open question of this subject and named it as the \emph{Riemann vs. MZ~derivatives conjecture.}

All these papers cite the current article as their primary source for the ideas originated in this work. And since the papers are a part of the new subject started in here, these pioneering ideas are more important now than they were when the article was completed.
As a clarification, the first version of this paper had the result of Corollary~\ref{C1} and its proof in a textual form. Based on this, subsequent papers such as \cite{AC2} and \cite{ACF1} have taken the corollary for granted and started directly with the statements of the variant and the GGR conjecture in their equivalent form given in Theorem~\ref{T2}(iii) and its similar statement for backward shifts instead of forward shifts. 

Up to the same publication delays described above, the present article is the fourth in a new wave of research aimed at solving major problems in generalized differentiation by using novel ideas and methods of proof coming from other areas of mathematics. Article~\cite{ACCs} that solved the main problem on the equivalence between Peano and single generalized Riemann derivatives, which we described in Section~0.1, used the theory of infinite linear systems of equations. Articles~\cite{ACCh} and~\cite{ACCH} solved the problem of finding the explicit description of all pairs $(\mathcal{A},\mathcal{B})$ of generalized Riemann differentiations for which $\mathcal{A}$-differentiation is equivalent to $\mathcal{B}$-differentiation, for all real or complex functions at a point, by using group algebras, group gradings, and symmetric functions. And this article uses a combinatorial algorithm in the proof of Theorem~\ref{T2} and an example based on group theory in the proof of Theorem~\ref{T0}(ii).

The same new wave of research is changing the focus on the equivalence between Peano and generalized Riemann derivatives, from an almost everywhere theory to a more precise pointwise theory. This brings more emphasis on the exact expressions of generalized Riemann derivatives and the relations between such expressions that translate into common properties of the derivatives. This is how algebraic and combinatorial methods found their way into the subject. The algebraic flavor started in \cite{ACCs} and continued in \cite{ACCh,ACCH,AC2,ACF1}. The combinatorial flavor started here and continued in \cite{AC2,ACF}.

\subsection*{Acknowledgment} Stefan Catoiu's research was supported in part by a one-quarter faculty research leave from the University Research Council at DePaul University during the fall of 2020.


\bibliographystyle{plain}

\subsection*{Statements and Declarations} The authors declare that there are no conflicts associated with this submission, and there is no data collection involved in the article.

\end{document}